\documentclass[10pt]{article}

\usepackage{cite}

\usepackage{amsfonts, amssymb, amsmath, amsthm, epsf}

\usepackage{graphicx}

\usepackage{fullpage}

\title{Uniqueness of transverse solutions for reaction-diffusion equations with spatially distributed hysteresis}
\author{Pavel Gurevich
\footnote{Free Univerity of Berlin, Peoples' Friendship University
of Russia, email: gurevichp@gmail.com} , Sergey Tikhomirov
\footnote{Free Univerity of Berlin, email:
sergey.tikhomirov@gmail.com} }

\theoremstyle{plain}
\newtheorem{theorem}{Theorem}[section]
\newtheorem{lemma}{Lemma}[section]

\newtheorem{condition}{Condition}[section]

\theoremstyle{definition}
\newtheorem{definition}{Definition}[section]

\newtheorem{remark}{Remark}[section]

\numberwithin{equation}{section} \numberwithin{figure}{section}

\renewcommand{\phi}{{\varphi}}

\newcommand{\cH}{{\mathcal H}}

\newcommand{\oQ}{{\overline Q}}
\newcommand{\oa}{{\overline a}}

\newcommand{\ox}{{\overline x}}
\newcommand{\ophi}{{\overline\varphi}}

\newcommand{\bbR}{{\mathbb R}}

\let\phi=\varphi

\begin{document}

\maketitle

\begin{abstract}
The paper deals with reaction-diffusion equations involving a
hysteretic discontinuity in the source term, which is defined at
each spatial point. Such problems describe biological processes
and chemical reactions in which diffusive and nondiffusive
substances interact according to hysteresis law. Under the
assumption that the initial data are spatially transverse, we prove
a theorem on the uniqueness of solutions. The theorem covers the
case of non-Lipschitz hysteresis branches arising in the theory of
slow-fast systems.
\end{abstract}



\section{Introduction}

We consider reaction-diffusion equations with  right-hand sides
involving a discontinuous hysteresis defined at each spatial
point. Such problems describe biological processes and chemical
reactions in which diffusive and nondiffusive substances interact
according to hysteresis law. As a result, various spatial and
spatio-temporal patterns may appear (see,
e.g.,~\cite{Jaeger1,Jaeger2}).

First rigorous results about the existence of solutions of
parabolic equations with hysteresis in the source term have been
obtained in~\cite{Alt, VisintinSpatHyst,Kopfova} for multi-valued
hysteresis. Formal asymptotic expansions of solutions were
recently obtained for some special case in~\cite{Ilin}. However,
the uniqueness of solutions and their continuous dependence on
initial data as well as a thorough analysis of pattern formation
remained open questions.

In~\cite{GurShamTikh}, a new approach was suggested. It allowed us
to find a broad class of initial data (transverse functions, see
Sec.~\ref{secSetting}) for which a solution exists and, if unique,
continuously depends on initial data. The approach is based on
tracking the so-called free boundary which defines the hysteresis
topology. The hysteresis topology is, in the simplest case,
related to the structure of subdomains in space where hysteresis
takes the same value. The main advantage of this approach is that,
tracking the hysteretic free boundary, one gets  necessary
information on the precise form of emerging spatio-temporal
patterns.

The transversality assumption roughly speaking means that the
initial function has a nonvanishing derivative on the boundary
between the above-mentioned subdomains at the initial moment (see
Condition~\ref{condA1'}). Under this assumption, we prove in the
present paper that the solution  for the reaction-diffusion
equation  with discontinuous spatially distributed hysteresis is
unique for some positive time. Combining this result with the
results in~\cite{GurShamTikh}, we obtain that the solution can be
uniquely continued as long as it remains transverse.

We consider a one-dimensional
domain and, for the clarity of exposition, we concentrate on a scalar reaction-diffusion equation. The paper is
organized as follows. In Sec.~\ref{secSetting}, we define
functional spaces, introduce spatially distributed hysteresis and
set the prototype problem. In the end of Sec.~\ref{secSetting}, we
formulate the main result of the paper: Theorem~\ref{tUniqueness}
on the uniqueness of transverse solutions.

Section~\ref{SecProofMainResult} is devoted to the proof of
Theorem~\ref{tUniqueness}.

Interestingly, the uniqueness of solutions holds for some classes
non-Lipschitz hysteresis branches, too. In particular, we prove
the uniqueness for hysteresis branches arising in bistable
slow-fast reaction-diffusion systems (see
e.g.,~\cite{Nishiura-SIMA-1990} and the references therein). In
Appendix~\ref{secSlowFast}, we briefly describe such systems and
show that the arising hysteresis branches satisfy
assumptions of our main theorem.

\section{Setting of the problem}\label{secSetting}

\subsection{Functional spaces}

We denote by $L_q=L_q(0,1)$, $q>1$, the standard Lebesgue space
and by $W_q^l=W_q^l(0,1)$ with natural $l$ the standard Sobolev
space. For a noninteger $l>0$, denote by $W_q^l=W_q^l(0,1)$ the
Sobolev space with the norm
$$
\|v\|_{W_q^l}=\|v\|_{W_q^{[l]}}+\left(\int_0^1dx\int_0^1
\dfrac{|v^{([l])}(x)-v^{([l])}(y)|^q}{|x-y|^{1+q(l-[l])}}dy\right)^{1/q},
$$
where $[l]$ is the integer part of $l$.

Let $Q_T=(0,1)\times (0,T)$. We introduce the H\"older space
$C^\gamma(\oQ_T)$, $0<\gamma<1$, and the anisotropic Sobolev space
$ W_q^{2,1}(Q_T) $ with the norm
$$
\|u\|_{W_q^{2,1}(Q_T)} =\left(\int_0^T \|u(\cdot,t)\|_{W_q^2}^q\,
dt+ \int_0^T \|u_t(\cdot,t)\|_{L_q}^q\, dt\right)^{1/q}.
$$

Throughout the paper, we fix $q$ and $\gamma$ such that
$$
q>3, \quad 0<\gamma< 1-3/q.
$$
Then $u,u_x\in C^\gamma(\overline Q_T)$ whenever $u\in
W_q^{2,1}(Q_T)$ (see Lemma 3.3 in~\cite[Chap. 2]{LadSolUral}).

In what follows, we will consider solutions of parabolic problems
in the space $W_q^{2,1}(Q_T)$. To define the space of initial
data, we will use the fact that if $u\in W_q^{2,1}(Q_T)$, then the
trace $u|_{t=t_0}$ is well defined and belongs to $W_q^{2-2/q}$
for all $t_0\in[0,T]$ (see Lemma 2.4 in~\cite[Chap.
2]{LadSolUral}). Moreover, one can define the space
$W_{q,N}^{2-2/q}$ as the subspace of functions from $W_q^{2-2/q}$
with the zero Neumann boundary conditions.

\subsection{Hysteresis}

In this section, we introduce a hysteresis operator defined for
functions of time variable $t$. Then we extend the definition to a
spatially distributed hysteresis acting on a space of functions of
time variable $t$ and space variable $x$.

We fix two numbers $\alpha$ and $\beta$ such that $\alpha<\beta$.
The numbers $\alpha$ and $\beta$ will play a role of thresholds
for the hysteresis operator. Next, we introduce continuous
functions ({\it hysteresis branches})
$$
H_1:(-\infty,\beta]\mapsto\bbR,\qquad H_2:[\alpha,\infty)
\mapsto\bbR.
$$

We assume throughout that the following condition holds.

\begin{condition}\label{condH}
There is a number $\sigma\in[0,1)$ such that, for any $U>0$, there
exists $M=M(U)>0$ with the properties
\begin{equation}\label{eqCondH1}
|H_1(u)-H_1(\hat u)|\le \dfrac{M}{(\beta-u)^\sigma+(\beta-\hat
u)^\sigma}|u-\hat u|,\quad \forall u,\hat u\in[-U,\beta),
\end{equation}
\begin{equation}\label{eqCondH2}
|H_2(u)-H_2(\hat u)|\le \dfrac{M}{(u-\alpha)^\sigma+(\hat
u-\alpha)^\sigma}|u-\hat u|,\quad \forall u,\hat u\in(\alpha,U].
\end{equation}
\end{condition}

\begin{remark}
\begin{enumerate}
\item Any locally Lipschitz continuous functions $H_1(u)$ and
$H_2(u)$ satisfy Condition~\ref{condH}. Moreover, this condition
covers the important case where the hysteresis branches $H_1(u)$
and $H_2(u)$ are the stable parts of the curve $g(u,v)=0$ in the
slow-fast system~\eqref{eq1'}. They are not Lipschitz continuous
near the points $u=\beta$ and $\alpha$, respectively. But they
still satisfy Condition~\ref{condH} (see
Appendix~\ref{secSlowFast} for details).

\item On the other hand, any $H_1(u)$ and $H_2(u)$ satisfying
Condition~\ref{condH} are locally H\"older continuous with
exponent $1-\sigma$ on $(-\infty,\beta]$ and $[\alpha,\infty)$,
respectively. Furthermore, $H_1(u)$ and $H_2(u)$ are locally Lipschitz on the open intervals
$(-\infty, \beta)$ and $(\alpha, \infty)$, respectively. Therefore, they satisfy the assumptions
in~\cite{GurShamTikh}, where existence of solutions and their
continuous dependence on initial data were proved.
\end{enumerate}
\end{remark}

We fix $T>0$ and denote by $C_r[0,T )$ the linear space of
functions which are continuous on the right in $[0,T )$. For any
$\zeta_0\in\{1,2\}$ ({\it initial configuration of hysteresis})
and $g\in C[0,T]$ ({\it input}), we introduce the {\it
configuration} function
$$
\zeta:\{1,2\}\times C[0,T]\to C_r[0,T),\quad
\zeta(t)=\zeta(\zeta_0,g)(t)
$$
 as follows. Let
$X_t=\{t'\in(0,t]:g(t')=\alpha\ \text{or } \beta\}$. Then
$$
\zeta(0)=\begin{cases}
1 & \text{if } g(0 )\le \alpha,\\
2 & \text{if } g(0)\ge \beta,\\
\zeta_0 & \text{if } g(0)\in(\alpha,\beta)\\
\end{cases}
$$
and for $t\in(0,T]$
$$
\zeta(t)=\begin{cases}
\zeta(0) & \text{if } X_t=\varnothing,\\
1 & \text{if } X_t\ne\varnothing\ \text{and } g(\max X_t)=\alpha, \\
2 & \text{if } X_t\ne\varnothing\ \text{and } g(\max X_t)=\beta.
\end{cases}
$$

\begin{figure}[ht]
        \centering
        \includegraphics[height=120pt]{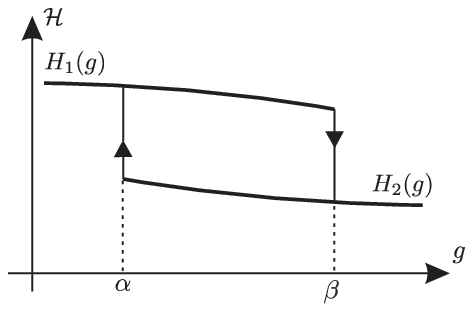}
        \caption{The hysteresis operator $\cH$}
        \label{figHyst}
\end{figure}

Now we introduce the {\it hysteresis operator\/}
(cf.~\cite{KrasnBook,Visintin})
$$
\cH: \{1,2\}\times C[0,T]\to C_r[0,T )
$$
by the following rule. For any initial configuration
$\zeta_0\in\{1,2\}$   and input $g\in C[0,T]$, the function
$\cH(\zeta_0,g):[0,T]\to\bbR$ ({\it output}) is given by
$$
\cH(\zeta_0,g)(t)=H_{\zeta(t)}(g(t)),
$$
where $\zeta(t)$ is the configuration function defined above (see
Fig.~\ref{figHyst}). Note that this hysteresis operator is typically not continuous.


Now we introduce a {\it spatially distributed hysteresis}. Assume
that the initial configuration and the input function depend on
spatial variable $x\in[0,1]$. Denote them by $\xi_0(x)$ and
$u(x,t)$, where
$$
\xi_0:[0,1]\mapsto \{1, 2\},\qquad u:[0,1]\times[0,T]\mapsto\bbR.
$$

Let $u(x,\cdot)\in C[0,T]$. Denote $\phi(x)=u(x,0)$. We say that a function $\varphi(x)$ and a hysteresis configuration $\xi_0(x)$ are
\textit{consistent} if, for any $x \in [0, 1]$,
$$ \xi_0(x)\in
\begin{cases}
\{1\} & \text{if } \phi(x)\ge\beta,\\
\{2\} & \text{if } \phi(x)\le\alpha,\\
\{1,2\} & \text{if } \phi(x)\in(\alpha,\beta).\\
\end{cases}
$$
Then we can define the function
\begin{equation}\label{eqSDH}
v(x,t)=\cH(\xi_0(x),u(x,\cdot))(t),
\end{equation}
which is called {\it spatially distributed hysteresis}.

\subsection{Reaction-diffusion equations with hysteresis}

The main object of this paper is the initial boundary-value
problem for the reaction-diffusion equation
\begin{align}
&u_t=u_{xx}+v,\qquad (x,t)\in Q_T, \label{eq1}\\
&u_x|_{x=0}=u_x|_{x=1}=0,\label{eq2}\\
&u|_{t=0}=\phi(x),\quad x\in(0,1),\label{eq3}
\end{align}
where $v=v(x,t)$ represents the spatially distributed hysteresis
given by~\eqref{eqSDH}.

\begin{remark}
The more general equation
$$u_t=u_{xx}+f(u,v)$$ with locally Lipschitz continuous right-hand
side $f$ can be reduced to Eq.~\eqref{eq1}. Indeed, it suffices to
replace the hysteresis branches $H_j(u)$ in the definition of
hysteresis $\cH$ by $ F_j(u)=f(u,H_j(u))$, $j=1,2$. One can check
that the functions $F_j(u)$ also satisfy Condition~\ref{condH}.
\end{remark}

The general assumption on the initial data $\phi(x)$ under which the
uniqueness result holds is that $\phi(x)$ is transverse with respect
to the initial configuration $\xi_0(x)$.

\begin{definition}\label{condTransverse}
We say that a function $\phi \in C^1[0, 1]$ is {\em transverse}
$($with respect to a spatial  configuration $\xi_0(x))$ if it is
consistent with $\xi_0(x)$ and the following holds:
\begin{enumerate}
\item if $\phi(\ox)=\alpha$ and $\phi'(\ox)=0$ for some
$\ox\in[0,1]$, then $\xi_0(\ox)=1$ in a neighborhood of $\ox$;

\item if $\phi(\ox)=\beta$ and $\phi'(\ox)=0$ for some
$\ox\in[0,1]$, then $\xi_0(\ox)=2$ in a neighborhood of $\ox$.
\end{enumerate}
\end{definition}

\begin{figure}[ht]
\begin{minipage}[b]{0.45\linewidth}
        \centering
        \includegraphics{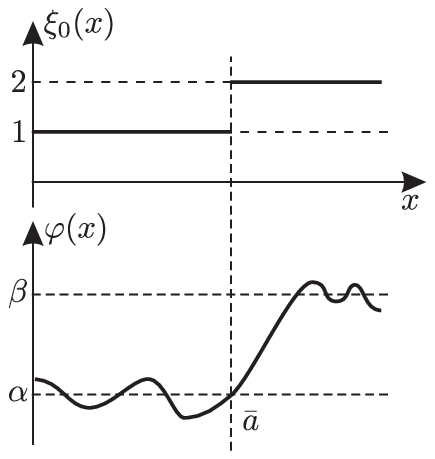}
        \caption{Initial data  satisfying Condition~\ref{condA1'}}
        \label{fig11}
\end{minipage}
\hspace{0.5cm}
\begin{minipage}[b]{0.45\linewidth}
        \centering
        \includegraphics{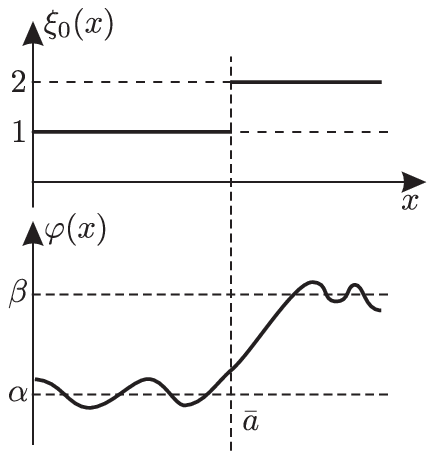}
        \caption{Initial data  satisfying the assumptions of Remark~\ref{rem*}}
        \label{fig11'}
\end{minipage}

\end{figure}

For the clarity of exposition, we now restrict ourselves to the
following prototype situation (see Fig.~\ref{fig11}). Fix some
$\oa\in (0,1)$.

\begin{condition}\label{condA1'}
\begin{enumerate}
\item For  $\oa\in(0,1)$, one has
$$
\xi_0(x)=\begin{cases}
1, & x\le\oa,\\
2, & x>\oa.
\end{cases}
$$
 \item The equation $\phi(x)=\beta$ on the interval $[0,\oa]$ has
 no roots.

\item
  The equation $\phi(x)=\alpha$ on the interval $[\oa,1]$ has
the unique root $x=\oa$.

\item $\phi'(\oa)>0$.

\end{enumerate}
\end{condition}

It follows from this condition that the hysteresis~\eqref{eqSDH} at
the initial moment is given by
$$
v|_{t=0}=\begin{cases}
H_1(\phi(x)), & x\le\oa,\\
H_2(\phi(x)), & x>\oa.
\end{cases}
$$

We give  a definition of a solution of
problem~\eqref{eq1}--\eqref{eq3},   assuming  that $\phi\in
W_{q,N}^{2-2/q}$.

\begin{definition}
A function $u(x,t)$ is called a {\em solution  of
problem~\eqref{eq1}--\eqref{eq3} {\rm (}in $Q_T${\rm) }} if $u\in
W_q^{2,1}(Q_T)$, $v(x,t)$ is measurable, $u$ and $v$ satisfy
equation~\eqref{eq1} for a.e. $(x,t)\in Q_T$, and
conditions~\eqref{eq2} and~\eqref{eq3} are satisfied in the sense
of traces.
\end{definition}

It follows from this definition that any solution $u(x,t)$ belongs
to $C(\oQ_T)$. Therefore, the function $v(x,t)$ is well defined
by~\eqref{eqSDH} and belongs to $L_\infty(Q_T)$.

\begin{remark}\label{rem*}
An initial function may also satisfy $\varphi(\bar{a}) > \alpha$.
Such a function is also consistent with $\xi_0(x)$ from Condition \ref{condA1'},
provided that the equation $\varphi(x) = \beta$ has no roots on the interval $[0, \bar{a}]$
and the equation $\varphi(x) = \alpha$ has no roots on the interval $[\bar{a}, 1]$ (see Fig. \ref{fig11'}).
However, this situation is much simpler.
Indeed, one can show (see \cite{GurShamTikh}) that, for any solution of problem \eqref{eq1}--\eqref{eq3}
with such an initial function, the hysteresis $\cH(\xi_0(x),u(x,\cdot))(t)$ does not depend on $t$ on a sufficiently
small time interval. Therefore, in this case, we leave details to the reader and rather concentrate on Condition \ref{condA1'}.
\end{remark}

\subsection{Main result}

In what follows, we always assume that Conditions~$\ref{condH}$
and~$\ref{condA1'}$ hold.

In~\cite{GurShamTikh}, it was proved that any solution of
problem~\eqref{eq1}--\eqref{eq3} remains transverse on some time
interval $[0,T]$ in the sense of the following definition.

\begin{definition}
A function $u\in C^{1,0}(\overline Q_T)$ is {\em
transverse} {\em on} $[0,T]$ $($with respect to a spatial
configuration $\xi(x,t))$ if, for every fixed $t\in[0,T]$, the
function $u(\cdot,t)$ is transverse with respect to the spatial
configuration $\xi(\cdot,t)$.
\end{definition}

Furthermore, problem~\eqref{eq1}--\eqref{eq3} has
at least one  transverse solution on the above interval $[0, T]$ and any such
solution can be extended to its maximal interval of existence on
which it remains transverse.

In this paper, we show that  there exists no more than one
transverse solution of problem~\eqref{eq1}--\eqref{eq3}. We
formulate the main result as follows.

\begin{theorem}\label{tUniqueness}
Let Conditions~$\ref{condH}$ and~$\ref{condA1'}$ hold, and let
$q>3$. Assume that $u,\hat u\in W_q^{2,1}(Q_{T_0})$ are two
transverse solutions of problem~\eqref{eq1}--\eqref{eq3} in
$Q_{T_0}$ for some $T_0$. Then $u=\hat u$.
\end{theorem}

\section{Proof of Theorem~\ref{tUniqueness}}\label{SecProofMainResult}

Taking into account Remark \ref{rem*}, it suffices to prove the uniqueness on a sufficiently
small time interval $[0,T]$, $T\le T_0$.

We denote
$$
\ophi=\dfrac{\phi_x(\oa)}{2}.
$$
By Condition~\ref{condA1'}, $\ophi>0$.

\begin{lemma}\label{lExistsab} Let $u\in W_q^{2,1}(Q_{T_0})$ be a
solution of problem~\eqref{eq1}--\eqref{eq3} for some $T_0>0$.
Then there exist $T\in(0,T_0)$ and $\delta>0$ such that the
following hold on the interval $t\in[0,T]$$:$
\begin{enumerate}
\item $u_x(x,t)\ge\ophi$ for all $x\in[\oa-\delta,\oa+\delta];$

\item the equation $u(x,t)=\alpha$ on the interval $x\in
[\oa-\delta,1]$ has a unique root $x=a(t);$

\item $a(t)$ is continuous$;$

\item the function
\begin{equation}\label{eqb0maxa0}
b(t)=\max\limits_{s\in[0,t]}a(s)
\end{equation}
satisfies $b(t)\in[\oa,\oa+\delta];$

\item the equation $u(x,t)=\beta$ on the interval $x\in[0,b(t)]$
has no roots.
\end{enumerate}
\end{lemma}
\proof The assertions of the lemma follow from the fact that
$u,u_x\in C^\gamma(\oQ_{T})$ for all sufficiently small $T$ (with
the norms in $C^\gamma(\oQ_{T})$ bounded uniformly with respect to
small $T$),
 $u(x,0)=\phi(x)$, and $\phi(x)$ is transverse (see
 Condition~\ref{condA1'}).
\endproof

By possibly decreasing $T$ and $\delta$, we see that
Lemma~\ref{lExistsab}  holds for $\hat u$ with some functions
$\hat a(t)$ and $\hat b(t)=\max\limits_{s\in[0,t]}\hat a(s)$
instead of $a(t)$ and $b(t)$.

Now the key observation is that the hysteresis $\cH$ acting on
transverse functions $u$ and $\hat u$ on the time interval $[0,T]$
can be represented in terms of the free boundaries $b(t)$ and
$\hat b(t)$ as follows:
\begin{equation}\label{eqHFreeBoundary}
\begin{aligned}
\cH(\xi_0(x),u(x,\cdot)(t)&=\begin{cases} H_1(u(x,t)), & 0\le x\le  b(t),\\
H_2(u(x,t)), & b(t)<x\le 1,
\end{cases}\\
\cH(\xi_0(x),\hat u(x,\cdot)(t)&=\begin{cases} H_1(\hat u(x,t)), & 0\le x\le  \hat b(t),\\
H_2(\hat u(x,t)), & \hat b(t)<x\le 1.
\end{cases}
\end{aligned}
\end{equation}

The following lemma allows us to estimate the distance between the
free boundaries $b(t)$ and $\hat b(t)$ in terms of the difference
between $u$ and $\hat u$.

\begin{lemma}\label{lb-hatb}
Let $T$ be the number from Lemma~$\ref{lExistsab}$. Then
$$
\|b-\hat b\|_{C[0,T]}\le \dfrac{1}{\ophi}\|u-\hat u\|_{C(\oQ_T)}.
$$
\end{lemma}
\proof It follows from the definition of $b(t)$ and $\hat b(t)$ that
$$
\|b-\hat b\|_{C[0,T]}\le \|a-\hat a\|_{C[0,T]}.
$$
On the other hand, using Lemma~\ref{lExistsab}, we obtain for any
$t\in[0,T]$
$$
|a(t)-\hat a(t)|\le \dfrac{1}{\ophi}|u(a(t),t)-\hat u(a(t),t)|\le
\dfrac{1}{\ophi}\|u-\hat u\|_{C(\oQ_T)},
$$
which completes the proof.
\endproof

Now we can prove Theorem~\ref{tUniqueness}.

1. Denote $w=u-\hat u$. The function $w$ satisfies the linear
parabolic equation
\begin{equation}\label{eqLinearw}
w_t=w_{xx}+h(x,t),\qquad (x,t)\in Q_T,
\end{equation}
where $h(x,t)=\cH(\xi_0(x),u(x,\cdot))-\cH(\xi_0(x),\hat
u(x,\cdot))$, and the zero boundary and initial conditions.
Obviously, $h\in L_\infty(Q_T)$, and the function $w$ can be
represented via the Green function $G(x,y,t,s)$ of the heat
equation with the Neumann boundary conditions:
\begin{equation*}
w(x,t)=\int\limits_0^t\int\limits_0^1 G(x,y,t,s)h(y,s)\,dy ds.
\end{equation*}
Therefore, using the estimate
$$
|G(x,y,t,s)|\le \dfrac{k_1}{\sqrt{t-s}}
e^{-{(x-y)^2}/{(4(t-s))}}\le \dfrac{k_1}{\sqrt{t-s}},\quad 0<s<t,
$$
with $k_1>0$  not depending on $(x,t)\in Q_T$ and $T>0$ (see,
e.g.,~\cite{Ivasisen}),
 we obtain
\begin{equation}\label{eqUniqueness1}
|w(x,t)|\le
k_1\int\limits_0^t\dfrac{ds}{\sqrt{t-s}}\int\limits_0^1
h(y,s)\,dy.
\end{equation}

2. Let us estimate the interior integral in~\eqref{eqUniqueness1}
for a fixed $s$. We assume that $b(s)<\hat b(s)$ (the case
$b(s)\ge\hat b(s)$ is treated analogously). Then, due
to~\eqref{eqHFreeBoundary},
$$
h(y,s)=
\begin{cases}
H_1(u)-H_1(\hat u),& 0<y< b(s),\\
H_2(u)-H_1(\hat u),& b(s)<y< \hat b(s),\\
H_2(u)-H_2(\hat u),& \hat b(s)<y< 1.
\end{cases}
$$

2.1. Assertion 5 in Lemma~\ref{lExistsab} implies that
$$
u(y,s)<\beta,\quad \hat u(y,s)<\beta
$$
on the closed set $\{(y,s): y\in[0, b(s)],\ s\in[0,T]\}$. Hence, the values $\beta-u(y,s)$ and $\beta-\hat
u(y,s)$ are separated from 0. Therefore, using Condition~\ref{condH}, we obtain

\begin{multline}\label{eqUniqueness2}
\int\limits_0^{b(s)} |h(y,s)|\,dy \le
\int\limits_{0}^{b(s)}\dfrac{M}{(\beta-u(y,s))^\sigma+(\beta-\hat
u(y,s))^\sigma}|u(y,s)-\hat u(y,s)|\,dy\\
 \le k_2 \int\limits_0^{b(s)}|u(y,s)-\hat u(y,s)|\,dy
\le k_2\|w\|_{C(\oQ_T)},
\end{multline}
where $k_2>0$ and the constants $k_3,k_4,k_5>0$ below do not
depend on $s\in[0,T]$.

2.2. By the boundedness of $H_1(\hat u)$ and $H_2(u)$ for
$(y,s)\in\oQ_T$, we have
\begin{equation*}
\int\limits_{b(s)}^{\hat b(s)} |h(y,s)|\,dy\le k_3
\int\limits_{b(s)}^{\hat b(s)}\,dy\le k_3\|b-\hat b\|_{C[0,T]}.
\end{equation*}
Applying Lemma~\ref{lb-hatb} yields
\begin{equation}\label{eqUniqueness3}
\int\limits_{b(s)}^{\hat b(s)} |h(y,s)|\,dy\le \dfrac{k_3}{\ophi}
\|w\|_{C(\oQ_T)}.
\end{equation}

\begin{figure}[ht]
        \centering
        \includegraphics[height = 40mm]{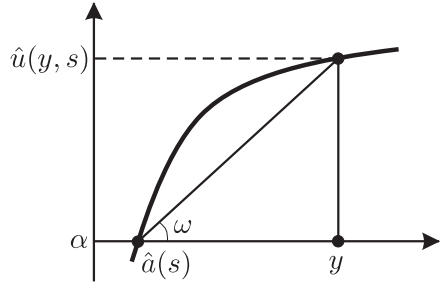}
        \caption{The mean-value theorem for $\hat u(y,s)$: $\tan \omega \ge\ophi$}
        \label{figualpha}
\end{figure}

2.3. Let $\delta$ be the number from Lemma~\ref{lExistsab}, and
let $\hat b(s)<y<\oa+\delta$. Then, using assertion~1 in
Lemma~\ref{lExistsab} and the mean-value theorem, we have (see
Fig.~\ref{figualpha})
$$
|\hat u(y,s)-\alpha|=\hat u(y,s)-\hat u(\hat a(s),s)\ge (y-\hat
a(s))\ophi\ge (y-\hat b(s))\ophi.
$$
Similarly,
$$
|u(y,s)-\alpha| \ge (y-\hat b(s))\ophi.
$$
Taking into account these two inequalities and using
Condition~\ref{condH}, we obtain

\begin{equation}\label{eqUniqueness4}
\begin{aligned}
\int\limits_{\hat b(s)}^{\oa+\delta} |h(y,s)|\,dy &\le
\dfrac{M}{2} \int\limits_{\hat b(s)}^{\oa+\delta}
\dfrac{|u(y,s)-\hat u(y,s)|}{(y-\hat b(s))^\sigma}\,dy\\
 &\le
\dfrac{M \|w\|_{C(\oQ_T)}}{2} \int\limits_{\hat b(s)}^{\oa+\delta}
\dfrac{1}{(y-\hat b(s))^\sigma}\,dy 
\le k_4\|w\|_{C(\oQ_T)}.
\end{aligned}
\end{equation}

2.4. Finally, assertions 1 and 2 in Lemma~\ref{lExistsab} imply
that
$$
u(y,s)>\alpha,\quad \hat u(y,s)>\alpha
$$
on the closed set $[\oa+\delta,1]\times [0,T]$. Hence, the values $u(y,s) - \alpha$ and $\hat u(y,s) - \alpha$ are separated from 0. Therefore, due to
Condition~\ref{condH},
\begin{equation}\label{eqUniqueness5}
\begin{aligned}
\int\limits_{\oa+\delta}^1 |h(y,s)|\,dy & \le
\int\limits_{\oa+\delta}^1\dfrac{M}{(u(y,s)-\alpha)^\sigma+(\hat
u(y,s)-\alpha)^\sigma}|u(y,s)-\hat u(y,s)|\,dy\\
&\le k_5 \int\limits_{\oa+\delta}^1|u(y,s)-\hat u(y,s)|dy\le
k_5\|w\|_{C(\oQ_T)}.
\end{aligned}
\end{equation}

3. Combining estimate~\eqref{eqUniqueness1} with
inequalities~\eqref{eqUniqueness2}--\eqref{eqUniqueness5}, we
obtain
$$
|w(x,t)|\le k_6 \|w\|_{C(\oQ_T)} \int_0^t \dfrac{ds}{\sqrt{t-s}}=2
k_6 T^{1/2}\|w\|_{C(\oQ_T)}.
$$
Since $k_6$ does not depend on $x$ and $t$, the latter inequality
is equivalent to $\|w\|_{C(\oQ_T)}\le 2k_6
T^{1/2}\|w\|_{C(\oQ_T)}$. Since $k_6$ does not depend on (small)
$T$ either, it follows that $w=0$ provided that $T>0$ is small
enough. \qed

\appendix

\section{Connection with slow-fast systems}\label{secSlowFast}

\begin{figure}[ht]
        \centering
        \includegraphics{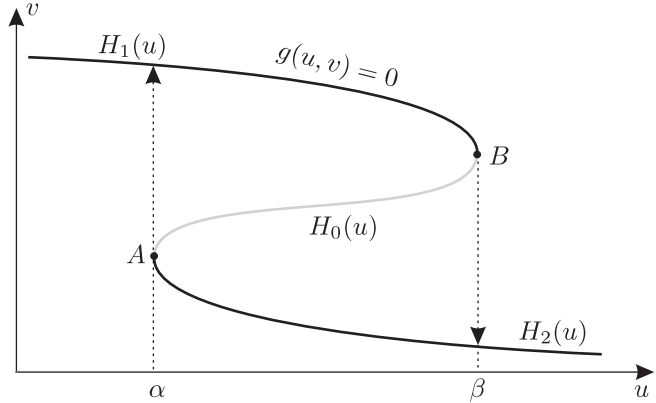}
        \caption{The nullcline of $g(u,v)$}
        \label{figuvfg}
\end{figure}

It is known~\cite{KrejciSlowFast, MischRozov} that
  hysteresis   may approximate solutions of  ordinary differential
equations with a small parameter $\varepsilon>0$ and a
``bistable'' right-hand side. Combining such an ordinary
differential equation with the reaction-diffusion
equation~\eqref{eq1}, one obtains a slow-fast system of the
 form
\begin{equation}\label{eq1'}
\left\{\begin{aligned} & u_t=u_{xx}+f(u,v),\\
& \varepsilon v_t=g(u,v),
\end{aligned}\right.
\end{equation}
where $\varepsilon>0$, $f(u,v)$ and $g(u,v)$ are smooth functions of class $C^{\infty}$,
and $g(u,v)$ satisfies the following condition (see
Fig.~\ref{figuvfg}).

\begin{condition}\label{condSlowFastg}
\begin{enumerate}
\item There are two numbers $\alpha<\beta$ such that the equation
$g(u,v)=0$ has a unique root $v=H_1(u)$ for $u<\alpha$ and
$v=H_2(u)$ for $u>\beta;$ three distinct roots $v=H_1(u)$,
$v=H_2(u)$, and $v=H_0(u)$ for $u\in(\alpha,\beta);$ two distinct
roots $v=H_1(\alpha)$ and $v=H_2(\alpha)=H_0(\alpha);$ two
distinct roots $v=H_1(\beta)=H_0(\beta)$ and $v=H_2(\beta)$.

\item the functions $H_1(u)$ and $H_2(u)$ are locally Lipschitz
continuous for $u<\beta$ and $u>\alpha$, respectively.

\item Let $A=(\alpha,H_2(\alpha))$ and $B=(\beta,H_1(\beta))$.
Then $\dfrac{\partial g(A)}{\partial u}\ne 0$ and $\dfrac{\partial
g(B)}{\partial u}\ne 0$.

\item $\dfrac{\partial g(A)}{\partial
v}=\dots=\dfrac{\partial^{n-1}g(A)}{\partial v^{n-1}}=0$,
$\dfrac{\partial^n g(A)}{\partial v^n}\ne 0$ for some even $n\ge
2$ and the same holds at the point $B$.

\item $g(u,v)>0$ $(<0)$ if a point $(u,v)$ lies to the left
$($right$)$ from the curve $g(u,v)=0$ on the plane $(u,v)$.
\end{enumerate}
\end{condition}

The simplest example of such a bistable nonlinearity typical,
e.g., for the FitzHugh--Nagumo systems  is given by
$$
g(u,v)=u+v-v^3.
$$

As we have mentioned before, the natural conjecture is that the
solutions of the slow-fast system~\eqref{eq1'}   approximate as
$\varepsilon\to 0$ the solution of Eq.~\eqref{eq1} with hysteresis
$\cH$ defined by the curves $H_1(u)$ and $H_2(u)$ from
Condition~\ref{condSlowFastg}. We refer
to~\cite{Nishiura-SIMA-1990}, where singular limit analysis of
traveling waves in bistable systems is done, and
to~\cite{Plotnikov, Evans}, where equations of the form
$u_t=\Delta\Phi(u)$ with a cubic nonlinearity $\Phi$ are studied.

The proof of the above conjecture in the general situation is an
open question, which is beyond the scope of the present paper.
However, it is clear that the proof would deal with hysteresis
defined by non-Lipschitz curves. The following result assures that
the non-Lipschitz curves $H_1(u)$ and $H_2(u)$ from
Condition~\ref{condSlowFastg}   satisfy Condition~\ref{condH} and
thus fit into our theory.

\begin{lemma}\label{lgHG}
Let $g(u,v)$ satisfy Condition~$\ref{condSlowFastg}$. Then the
functions $H_1(u)$ and $H_2(u)$ from
Condition~$\ref{condSlowFastg}$ satisfy Condition~$\ref{condH}$
with $\sigma=(n-1)/n$.
\end{lemma}
\proof 1. Without loss of generality, we assume that $A=0$ and
prove the lemma for the function $H_2(u)$. Due to item 2 in
Condition~\ref{condSlowFastg}, it suffices to prove
inequality~\eqref{eqCondH2} for $u$ in a small neighborhood of
$\alpha=0$.

Since $\dfrac{\partial g(A)}{\partial u}\ne 0$, it follows from
the implicit function theorem that there is a unique function
$G_2(v)$ defined in a small neighborhood of $0$ such that
$g(G_2(v),v)\equiv 0$ and
\begin{equation}\label{lgHG1}
G_2(0)=G_2'(0)=\dots =G_2^{(n-1)}(0)=0,\quad G_2^{(n)}(0)>0.
\end{equation}
Obviously, $G_2(v)$ and $H_2(u)$ are inverse to each other for all
negative $v$ and positive $u$ close to $0$. It is convenient to
introduce the function $G(w)=G_2(-w)$ defined for small positive
$w$. Due to~\eqref{lgHG1}, it satisfies
\begin{equation}\label{lgHG2}
G(0)=G'(0)=\dots G^{(n-1)}(0)=0,\quad G^{(n)}(0)>0
\end{equation}
since $n$ is even. Now if we denote $w=-H_2(u)$ ($>0$), then
$u=G_2(H_2(u))=G_2(-w)=G(w)$. Therefore, the inequality for $H_2$
in Condition~\eqref{condH}, which we have to prove, is equivalent
to the following:
\begin{equation}\label{lgHG3}
\dfrac{G(w)-G(\hat w)}{w-\hat w}\ge M\left
([G(w)]^\frac{n-1}{n}+[G(\hat w)]^\frac{n-1}{n}\right)
\end{equation}
for all $0<\hat w< w\le \varepsilon_0$, where $\varepsilon_0>0$ is
sufficiently small and $M>0$ does not depend on $w,\hat w$.

2. Expanding $G(w)$ by the Taylor formula about $w=\hat w$, we
have
\begin{equation}\label{lgHG4}
\dfrac{G(w)-G(\hat w)}{w-\hat w}=G'(\hat
w)+\sum\limits_{k=2}^{n-1} \dfrac{G^{(k)}(\hat w)}{k!}(w-\hat
w)^{k-1}+\dfrac{G^{(n)}(\xi)}{n!}(w-\hat w)^{n-1},
\end{equation}
where $0\le\xi=\xi(w,\hat w)\le\varepsilon_0$.

Further, we expand $G^{(k)}(\hat w)$, $k=1,\dots,n-1$ by the
Taylor formula about the origin, using~\eqref{lgHG2}:
$$
G^{(k)}(\hat w)=\dfrac{G^{(n)}(\xi_k)}{(n-k)!}\hat w^{n-k},
$$
where $0\le\xi_k=\xi_k(\hat w)\le\varepsilon_0$. Substituting
these expressions into~\eqref{lgHG4} yields
\begin{equation}\label{lgHG5}
\begin{aligned}
\dfrac{G(w)-G(\hat w)}{w-\hat w}
&\ge\dfrac{G^{(n)}(\xi_1)}{(n-1)!}\hat w^{n-1}
+\dfrac{G^{(n)}(\xi)}{n!}(w-\hat
w)^{n-1}\\
&\ge  \dfrac{G^{(n)}(0)}{2(n-1)!}\left(\hat
w^{n-1}+\dfrac{1}{n}(w-\hat w)^{n-1}\right),
\end{aligned}
\end{equation}
where we assume $\varepsilon_0$ so small that $G^{(n)}(\xi_k)\ge
\dfrac{G^{(n)}(0)}{2}$ and $G^{(n)}(\xi)\ge
\dfrac{G^{(n)}(0)}{2}$.

It is easy to show that 
$$
\hat{w}^{n-1} + \frac{1}{n}(w - \hat{w})^{n-1} \geq M_1 \left(w^{n-1} + \hat{w}^{n-1}\right)
$$
for some $M_1 > 0$ not depending on $w$ and $\hat{w}$, $0<\hat{w}<w$. Hence, inequality~\eqref{lgHG5} implies
\begin{equation}\label{lgHG6}
\begin{aligned}
&\dfrac{G(w)-G(\hat w)}{w-\hat w}\ge M_2\left(w^{n-1}+\hat
w^{n-1}\right),
\end{aligned}
\end{equation}
where $M_2>0$ does not depend on $w$ and $\hat{w}$, $0< \hat{w} < w \le \varepsilon$.
Combining~\eqref{lgHG6} and the relations
$$
G(w)=\dfrac{G^{(n)}(\zeta)}{n!}w^n\le \dfrac{2G^{(n)}(0)}{n!}w^n
$$
 with $0\le\zeta=\zeta(w)\le\varepsilon_0$
yields~\eqref{lgHG3}.
\endproof

\bigskip

{\bf Acknowledgement:} The authors are grateful to Willi J\"ager
for drawing their attention to the field of hysteresis.
The research of the first author was supported by the DFG project
SFB 910, by the DAAD program G-RISC, and by the RFBR (project
10-01-00395-a). The research of the second author was supported by
the Alexander von Humboldt Foundation.

\end{document}